\magnification=\magstep1

{\bf Eine Bemerkung zu einem Satz von E. Becker und D. Gondard}

\vskip5mm

J.-L. Colliot-Th\'el\`ene

\vskip5mm

{\it Eberhard Becker zum 60. Geburtstag gewidmet}

\vskip5mm

In ihrer Arbeit [B-G] geben E. Becker und D. Gondard
  eine algebraische Formel f\"ur die Anzahl von
Zusammenhangskomponenten des Raumes der reellen Punkte
einer reellen, glatten, projektiven Variet\"at.
In der Arbeit [C-T] hatte ich eine andere
algebraische Formel f\"ur diese Anzahl
gegeben.
In dieser Note zeige ich, wie man von einer Formel
zur anderen gehen kann  -- ohne von diesen beiden S\"atzen
Gebrauch zu machen.

  Becker und Gondard benutzen den Raum aller
${\bf R}$-Stellen eines K\"orpers, sowie einen
Satz von Ludwig Br\"ocker diesen Raum betreffend.
Dieser Raum wird
hier nicht benutzt. Daf\"ur wird
ein Reinheitssatz von Markus Rost angewandt.
Sowohl in [B-G] als in dieser Note wird Beckers
Charakterisierung der Summen $n$-ter Potenzen in einem K\"orper
benutzt.

\bigskip

  Zitieren wir zuerst die erw\"ahnten S\"atze.

Sei $X/{\bf R}$ eine glatte, projecktive, absolut
irreduzible Variet\"at \"uber dem K\"orper ${\bf R}$
der (\"ublichen) reellen Zahlen.
Sei $X({\bf R})$ die Menge der reellen Punkte von $X$.
Nehmen wir an, da{\ss} diese Menge nicht leer ist.
Sei $S$ die (endliche) Menge der
Zusammenhangskomponenten
von $X({\bf R})$ und $s$ die Ordnung von $S$.

Sei $K={\bf R}(X)$ der Funktionenk\"orper von $X$.
Sei $D(X) \subset K^*$ die multiplikative
Untergruppe aller Funktionen in ${\bf R}(X)^*$, die sich  in jedem
Punkt
$P$ von
$X$ (im schematischen Sinne)
als Produkt einer Einheit im lokalen Ring $O_{X,M}$
und einer Summe von
Quadraten von Elementen von $K$ schreiben lassen.

\bigskip

{\bf Satz 1} ([CT]) \hskip2mm {\it Der Quotient
$D(X)/(K^* \cap (\sum K^2)) $
  ist isomorph zur Gruppe $({\bf
Z}/2)^S$.}

\bigskip

{\bf Satz 2 } (Becker und Gondard, [B-G]) \hskip2mm {\it  Der Quotient
$(K^{*2} \cap (\sum K^4))/(K^* \cap (\sum K^2)^2) $  ist endlich, von
der
Ordnung $2^{s-1}$. }

\medskip

Da{\ss} $K^* \cap (\sum K^2)^2$ eine Untergruppe von
$K^{*2} \cap (\sum K^4) $ ist, folgt aus dem
Satz :

\medskip

{\bf Satz 3} (Becker, [B]) \hskip2mm {\it Sei $K$ ein K\"orper. Sei
$n$
eine gerade Zahl.
Ein Element $f \in K^*$ ist eine Summe von $n$-ter Potenzen
in $K$ dann und nur dann, wenn die  folgenden beiden
Bedingungen  erf\"ullt sind:

(a) Das Element $f$ ist eine Summe von Quadraten in $K$.

(b) F\"ur jede Krullbewertung $v$ von $K$ mit formalreellen
Restklassenk\"orper ist $v(f)$ durch $n$ teilbar.}

\bigskip

Die Quadratabbildung  $x \mapsto x^2$
induziert eine surjektive Abbildung
$$ K^*/K^* \cap (\sum K^2) \to K^{*2}/(K^{*2} \cap (\sum K^2)^2)    $$
mit Kern $\{\pm 1\}$.

Nach dem Satz von Becker hat man
$$ (\sum K^2)^2 \subset \sum K^4,$$
also
$$  K^* \cap(\sum K^2)^2 \subset K^{*2} \cap \sum K^4.$$

\medskip

{\bf Hauptsatz} \hskip2mm {\it Ein Element $f \in K^*$ liegt in der
Gruppe
$D(X)$ dann und nur dann, wenn das Element $f^2$  in $K^{*2} \cap \sum
K^4$ liegt.}

\medskip

Aus diesem Satz folgt sofort, da{\ss} die Abbildung $x \mapsto x^2$ 
eine
exakte Folge induziert
$$ 1 \to \{\pm 1\} \to D(X)/(K^* \cap \sum K^2)
\to (K^{*2} \cap \sum K^4)/(K^* \cap \sum K^2)^2 \to 1,$$
  die  die Verbindung zwischen Satz 1 und Satz 2 genau erkl\"art.

\bigskip

{\it Beweis des Hauptsatzes}. Sei $f \in D(X)$. Sei $v$ eine
Krullbewertung  von $K$ mit formalreellen
Restklassenk\"orper. Sei $A$  der Bewertungsring.
Da $X/{\bf R}$ projektiv ist, besitzt $A$ ein Zentrum auf
$X$, das hei{\ss}t, es gibt einen (nicht unbedingt abgeschlossenen)
Punkt
  $M \in X$, so da{\ss} die nat\"urliche
Abbildung ${\rm Spec} \hskip0,5mm {\bf R}(X) \to X$ einen
Homomorphismus von lokalen Ringen $O_{X,M} \to A $ induziert.
Da $f$ in $D(X)$ liegt, kann man $f=u.g$ schreiben, mit
$u \in O_{X,M} ^*$ und $g \in \sum K^2$, also
  $f=u.g$ mit $u \in A^*$ und $g \in K^* \cap  \sum K^2$.
Da der Restklassenk\"orper von $A$ formalreell ist,
folgt $2 \mid v(f)$. Aber dann hat man $f^2 \in K^{*2}$
und $4 \mid v(f^2)$. Dies gilt f\"ur eine beliebige Krullbewertung $v$
von
$K$ mit formallreellem Restklassenk\"orper. Aus dem Satz von Becker
folgt $f^2 \in \sum K^4$.

\medskip

Sei umgekehrt $f \in K^*$ gegeben, mit der Eigenschaft
$f^2 \in K^{*2} \cap \sum K^4$.

Sei $v$ eine
  Krullbewertung von
$K$ mit formalreellem Restklassenk\"orper.
Dann hat man $4 \mid v(f^2)$, also $2 \mid v(f)$.
Wenn $M \in X$ ein Punkt der Kodimension 1
  mit formalreellem Restklassenk\"orper ist,
dann kann man $f$ als Produkt  einer Einheit in $O_{X,M}$
und eines Quadrates in $K$ schreiben.

Sei $M \in X$
ein Punkt der Kodimension 1
  mit nichtformalreellem Restklassenk\"orper.
Sei $\pi $ ein Erzeuger des maximalen Ideals des
diskreten  Bewertungsrings $R=O_{X,M}$.
In $R$ kann man Elemente $a_i, i=1,\dots,m$ und $b$ finden,
mit $1+\sum_i a_i^2=\pi. b \in R$. Aus dieser Gleichung folgt
die Gleichung $(1+\pi )^2+ \sum_i a_i^2 =\pi.b + 2\pi +\pi ^2  \in
R$. Wenn $b \in R $ keine Einheit ist, dann ist $2+b+\pi $ eine
Einheit.
Also kann man $\pi $ als Produkt  einer Einheit in $O_{X,M}$
und einer Summe von Quadraten in $K$ schreiben. Daraus folgt,
da{\ss} sich jedes Element in $K^*$ als Produkt  einer Einheit
in $O_{X,M}$ und einer Summe von Quadraten in $K$ schreiben l\"asst.

Also: wenn $f \in K^*$  die Eigenschaft hat, da{\ss} sein Quadrat
$f^2$ in $ K^{*2} \cap \sum K^4$ liegt, dann kann man $f$ in jedem
Punkt der Kodimension 1 als   Produkt  einer Einheit
in $O_{X,M}$ und einer Summe von Quadraten in $K$ schreiben.

Nach einem bekannten Satz von A. Pfister ([P]) ist jede Summe von
Quadraten
  in $K={\bf R}(X)$  eine Summe von $2^d$ Quadraten, wobei $d$
die Dimension von $X$ ist.

Jetzt kann man einen allgemeinen Reinheitssatz von Rost [R] in dieser
speziellen Situation anwenden: wenn $f \in {\bf R}(X)^*$ sich in
jedem Punkt
$M$ der Kodimension 1 der glatten Variet\"at $X$ als Produkt einer
Einheit in
$M$ und einer Summe von
$2^d$ Quadraten im Funktionenk\"orper von $X$ schreiben l\"asst, dann
besitzt
$f$ auch eine solche Darstellung  in jedem Punkt von $X$. Also
geh\"ort
$f$ der Gruppe $D(X)$ an.

\vskip1cm

{\it Literatur}

\medskip

[B] E. Becker, Summen $n$-ter Potenzen in K\"orpern, J. reine angew.
Math. (Crelle) {\bf 307/308} (1979) 8--30.

\smallskip

[B-G] E. Becker und D. Gondard, Notes on the space of real places
of a formally real field, in {\it Real analytic and
algebraic geometry} (Trento, 1992),  21--46, de Gruyter, Berlin, 1995.

(Einen kurzen Beweis f\"ur den Satz von Becker und Gondard
hat Claus Scheiderer gegeben :
http://www.uni-duisburg.de/FB11/FGS/F1/claus.html$^{\#}$notes)
\smallskip
[CT] J.-L. Colliot-Th\'el\`ene, Formes multiplicatives et vari\'et\'es
alg\'ebriques, Bull. Soc. math. France {\bf 106} (1978) 113--151.
\smallskip
[P] A. Pfister,
Zur Darstellung definiter Funktionen als Summe von Quadraten.
Invent. math. {\bf 4} (1967) 229--237.
\smallskip
[R] M. Rost, Durch Normengruppen  definierte
birationale Invarianten, C. R. Acad. Sci. Paris S\'er. I,
Math\'ematiques,  {\bf 310}  (1990),  no. 4, 189--192.

\vskip2cm

J.-L. Colliot-Th\'el\`ene

C.N.R.S.,
Math\'ematiques

UMR 8628

  B\^atiment 425

Universit\'e Paris-Sud

F-91405 Orsay

Frankreich

\smallskip

colliot@math.u-psud.fr 

\vskip2cm

\bye